\documentclass[journal]{IEEEtran}
\usepackage{cite}    
\usepackage{comment}    
\usepackage{graphicx}  
\usepackage{caption}
\usepackage{subcaption}
\usepackage{caption}

\usepackage[table,dvipsnames]{xcolor}  
\usepackage{amssymb}
\usepackage{mathrsfs}
\usepackage{mathtools}   

\usepackage[normalem]{ulem}
\usepackage{cancel}

\usepackage{hhline}
\usepackage{tabularx,colortbl}

\usepackage[sans]{dsfont}
\usepackage{stmaryrd}
\usepackage{pifont}
\usepackage{wasysym}
\usepackage{algorithm} 
\usepackage{array}
\usepackage{url}

\usepackage{tabu}
\usepackage{comment}
\usepackage{pifont}

\usepackage{multirow}
\usepackage{booktabs}
\usepackage{algpseudocode}
\usepackage{algorithm}

\newtheorem{theorem}{Theorem}[section]

\newtheorem{definition}[theorem]{Definition}
\newtheorem{example}{Example}

\usepackage{soul}

\newcommand{\cmark}{\ding{51}}%
\newcommand{\xmark}{\ding{55}}%

\newcommand{\re}{\ensuremath{\mathrm{Re}}}
\newcommand{\im}{\ensuremath{\mathrm{Im}}}
\def\CC{\mathbb C}
\def\RR{\mathbb R}

\def\cG{\mathcal G}

\def\cN{\mathcal N}

\newcommand{\hide}[1]{}
\newcommand{\raf}[1]{(\ref{#1})}

\newcommand{\cV}{\ensuremath{\mathcal{V}}}
\newcommand{\cE}{\ensuremath{\mathcal{E}} }

\newcommand{\Opt}{\ensuremath{\textsc{Opt}}}

\DeclareCaptionType{copyrightbox}

\newtheorem{innercustomthm}{\textbf{Theorem}}
\newenvironment{customthm}[1]
  {\renewcommand\theinnercustomthm{#1}\innercustomthm}
  {\endinnercustomthm}

  \newtheorem{innercustomlem}{\textbf{Lemma}}
\newenvironment{customlem}[1]
  {\renewcommand\theinnercustomlem{#1}\innercustomlem}
  {\endinnercustomthm}

\title{Efficient Algorithm for Scalable Event-based Demand Response Management in Microgrids}

\author{
Areg Karapetyan, Majid Khonji, Chi-Kin Chau {\em Member}, {\em IEEE}, Khaled Elbassioni, and\\ H. H. Zeineldin, {\em Senior Member}, {\em IEEE}
\thanks{A. Karapetyan, C.-K. Chau, and K. Elbassioni are with the Department of Electrical Engineering and Computer Science at Masdar Institute of Science and Technology, UAE (e-mail: \{akarapetyan, ckchau, kelbassioni\}@masdar.ac.ae).
M. Khonji is with Dubai Electricity and Water Authority (DEWA), Dubai, UAE (e-mail: majid.khonji@dewa.gov.ae).
H. H. Zeineldin is with the Electrical Engineering and Computer Science Department, Masdar Institute of Science and Technology, Abu Dhabi, United Arab Emirates, and currently on leave from the Faculty of Engineering, Cairo University, Giza 12316, Egypt. (email: hzainaldin@masdar.ac.ae)}
}

\begin{document}

\maketitle

\begin{abstract}
	Demand response management has become one of the key enabling technologies for smart grids. Motivated by the increasing demand response incentives offered by service operators, more customers are subscribing to various demand response programs. However, with growing customer participation, the problem of determining the optimal loads to be curtailed in a microgrid during contingencies within a feasible time frame becomes computationally hard. This paper proposes an efficient approximation algorithm for event-based demand response management in microgrids. In event-based management, it is important to curtail loads as fast as possible to maintain the stability of a microgrid during the islanded mode in a scalable manner. A simple greedy approach is presented that can rapidly determine a close-to-optimal load curtailment scheme to maximize the aggregate customer utility in milliseconds for a large number of customers. This paper further derives a novel theoretical guarantee of the gap between the proposed efficient algorithm and the optimal solution (that may be computationally hard to obtain). The performance of algorithm is corroborated extensively by simulations with up to thousands of customers. For the sake of practicality, the proposed event-based demand response management algorithm is applied to a feeder from the Canadian benchmark distribution system. The simulation results demonstrate that the proposed approach efficiently optimizes microgrid operation during islanded mode while maintaining appropriate voltage levels and network constrains.

\end{abstract}

{\keywords Demand response management, microgrids, approximation algorithms, integer variables.}

\vspace*{10pt}
\section*{Nomenclature}

\addcontentsline{toc}{section}{Nomenclature}
\begin{IEEEdescription}[\IEEEusemathlabelsep\IEEEsetlabelwidth{{\sc UMDR}}]
\item[${\cal N}$] Set of customers
\item[$n$] Cardinality of customers' set
\item[$k$] Index of a customer
\item[$u_k$] Utility of the $k$-th customer if his power demand is not curtailed
\item[$x_k$] Binary decision variable if the $k$-th customer's power demand is retained (i.e., not curtailed)
\item[$X$] Set of customers whose power demand are retained (i.e., not curtailed) 
\item[$S_k$] Apparent power demand of the $k$-th customer, represented by a complex number
\item[$P_k$] Active power demand of the $k$-th customer, represented by a real number
\item[$Q_k$] Reactive power demand of the $k$-th customer, represented by a real number
\item[$C$] Total (apparent) generation power capacity
\item[$\theta$] Maximum difference between the phase angles of any pair of power demands.
\item[{\sc UMDR}] Utility Maximizing Demand Response Problem
\end{IEEEdescription}

\section{Introduction}

\PARstart{M}{icrogrids} (MGs) typically are medium-to-low voltage networks with integrated Distributed Generation (DG), capable of operating in grid connected or islanded mode. There is a high probability that an MG once initiated will be short of power, consequently resulting in significant voltage and frequency deviations, and leading to its instability. Designing an MG with the capability of operating in an islanded mode can enhance system reliability and power quality.

Demand response (DR) programs can be broadly classified into three classes: economic DR, emergency DR, and ancillary services DR. In particular, emergency DR \cite{Ref6} is utilized when there is insufficient supply of power to meet the available demand, especially for MGs. DR is a key feature for smart grids and can be used to alter loads during contingency conditions.  DR has been proven to have many benefits including the reduced price variations \cite{Ref1}, increased reliability \cite{Ref2}, enhanced congestion management \cite{Ref3} and security enhancement \cite{Ref4}. 

Sudden islanding of MGs can cause high imbalances between the local generation and demand and thus, DR management strategies are required to ensure the endurance of an MG during its autonomous operation \cite{Ref10}. The DR management schemes for islanded MGs have received increasing attention recently \cite{6496178}. Innovative DR management strategies for MGs can contribute to the improvement of stability especially during emergency conditions \cite{Ref10}. The emergency DR methods, proposed in \cite{Ref10,Ref12}, are based on local frequency measurement to decide the on/off states of loads. However, such methods do not take into account customer utility. In \cite{Abdelaziz2014}, a multi-stage centralized scheme for MGs is proposed in order to minimize the load shedding cost. The problem is formulated as a non-linear integer programming problem with power flow equations and voltage limits, and is solved using a numerical solver in a small scale. In \cite{6716076}, a DR scheme is proposed to minimize the total amount of load shedding in order to resolve the transformer overload. The aforementioned prior studies concerning centralized demand management for MGs have only considered systems with a small number of loads (e.g., less than $50$), and thus the optimization of operation for an MG within a short time frame during emergency conditions (within milliseconds) is plausible. 

The energy management systems for control and optimization of MGs can be broadly classified into two groups, distributed and centralized \cite{hatziargyriou}. In \cite{nikos4454000, 1626398, 5440980} a centralized MG control strategy is envisioned with the goal of effectively optimizing system operation and/or economic benefits. Usually, centralized energy management approaches are more efficient and, unlike the distributed ones, do not require computations to be performed in each node of the network as well as do not face synchronization problems. Distributed strategies, on the other hand, do not suffer from the scalability problem and are computationally tractable for large number of customers. With growing customer participation, centralized energy management approaches for MGs require large and extensive computational resources, and thus a distributed energy management approach is adopted in several prior studies. A distributed energy management strategy that minimizes the dissatisfaction of customers and power losses is proposed in \cite{6980137}. In \cite{6840288} a scalable DR scheme is proposed due to a large number of customers who need to be managed. The main objective is to maximize the aggregate utilities of the appliances in the households and minimize power losses.

This paper proposes a scalable centralized event-based DR management algorithm for MGs considering a large number of customers. The proposed algorithm can pave the foundations for more effective distributed strategies. The problem is formulated as an integer programming problem of determining which loads to be curtailed, as to maximize the cumulative customer utility subject to the maximum available MG capacity. This problem is known to be computationally hard (in particular,  it is \textit{strongly} {\sc NP-hard} to obtain an optimal solution). In this paper, an efficient greedy algorithm is proposed to achieve very short computational time (within a few milliseconds), in spite of a large number of customers. An important aspect of this approach is that it is capable of obtaining a close-to-optimal solution, with a theoretical guarantee of the gap to an optimal solution. The proposed approach is empirically evaluated in terms of computational time and quality of solution, compared with the conventional numerical solver (e.g., Gurobi). The observed computational time of the proposed algorithm significantly outperforms that of the conventional numerical solver by an order of $10^3$ times faster. With regards to the observed quality of solution of the proposed algorithm depicted in Fig. \ref{fig:obj} and Fig. \ref{fig:cr} in Section \ref{simresutssec}, the empirical results of various practical case studies recorded up to $10\%$ deviation from the optimal solution generated by the numerical solver, considering up to $1000$ customers.
	
Lastly, we investigate the feasibility of the load curtailment schemes produced by the proposed DR management algorithm on a 4-bus feeder from Canadian benchmark distribution system. The observed simulation results show that the solutions produced by the proposed algorithm are feasible and do not violate network constraints and voltage levels.

\section{Problem Formulation} \label{sec:model}

As mentioned earlier, most of the available literature considers a small number of customers. As indicated in \cite{6980137}, a centralized controller for energy management will likely be computationally demanding. One of the main challenges for MGs is to maintain stable operations especially during the islanded mode when there is a shortage of power. This paper proposes a scalable framework for implementing an event-based DR scheme that can optimally determine the loads to be curtailed subject to a capacity constraint taking into account significant customer participation. In this scheme, each customer submits his reactive and active power demand, and his utility that reflects how much he is willing to pay if his demand is satisfied. This setting assumes a customer submitting a single utility value that reflects the value of retaining the complete set of own demands. It is noteworthy to mention, that this scheme could be easily extended for practical purposes such that a customer is allowed to submit different utilities (independent of each other) for corresponding subsets of own demands. This enables customers to prioritize own demands based on their importance and value without altering the proposed formulation. A computationally efficient centralized MG controller with an objective of maximizing the total utility of the satisfied customers is proposed. The MG controller will be responsible for sending curtailment signals to the corresponding customers. To effectively control load curtailment, low-latency communication infrastructure between the controller and customers using separate power supply will be utilized in event-based DR management. Such a low-latency communication infrastructure is enabled by the standards of smart grid communication protocols \cite{smartgridbook}.

A Utility Maximizing Demand Response (UMDR) problem is formulated with an objective to maximize the total utility of satisfiable customers:

\begin{eqnarray}
\textsc{(UMDR)} & &  \displaystyle \max_{x_k \in \{0, 1 \}} \sum_{k\in{\cal N}}u_k x_k \label{DAP}\\
\text{subject to}& &  \displaystyle \Big|\sum_{k\in {\cal N}}S_k x_k\Big| \le C \, , \label{maxPA}
\end{eqnarray}
where ${\cal N}$ is the set of customers, $u_k$ is the utility of $k$-th customer, $S_k = P_k + {\bf i} Q_k \in\CC$ is the {\em complex-valued} apparent power demand for $k$-th customer, $C \in\RR_+$ is a real-valued generation capacity of total apparent MG power. Here, $x_k$ is a binary decision variable that takes value $1$ if and only if the $k$-th customer's power demand is retained. Note that the customers' utilities may be provided in advance, which will be taken into account during the islanding of MGs. It is worthwhile mentioning that when deploying the DR scheme practically time-varying customer utilities could be incorporated. Particularly, this scheme relies on a dynamic rather than a flat customer utility plan in which utilities could fluctuate based on time of day, so empowering customers' choice over utility values for peak and off-peak hours. Upon initiation of the proposed algorithm, the respective utility value is adopted according to the execution timing. In MGs, most power demand can be attributed to customers' demands, and hence, the effect of transmission loss is negligible.


While {\sc UMDR} problem as a {\it Quadratically Constrained Integer Programming} (QCIP) problem is not entirely new \cite{woeginger2000does}, its application to power systems was recently realized in \cite{Yu:2013}. There were only few results in the literature concerning {\sc UMDR} problem unlike its well-known relaxations, namely the classical 0-1 knapsack and continuous knapsack problems. The former problem resembles the special case of {\sc UMDR} problem. Specifically, {\sc UMDR} problem is equivalent to the classical knapsack problem when setting zero reactive power, namely, $Q_k = 0$ for all $k \in {\cal N}$. This seemingly small change in the problem formulation, in effect, drastically impacts its computational complexity. Note that {\sc UMDR} is a quadratically constrained integer programming problem, whereas the classical knapsack is a \textit{linear integer programming problem}. The latter problem is known to be \textit{weakly} {\sc NP-hard}\cite{KPP10book}, in other words, it is highly unlikely that there is any fully polynomial time algorithm to compute the exact optimal solution.

On the other hand, the presence of complex-valued power demands in {\sc UMDR} problem creates substantially more challenging problem than the classical knapsack problem (which is shown to be \textit{strongly} {\sc NP-hard} \cite{CKM14,CKM2015,KCE2014Hard}). Though there exist in the literature a number of exact and approximation algorithms for the classical knapsack problem including greedy algorithms and dynamic programming \cite{KPP10book}, there is no known efficient practical algorithm that can solve {\sc UMDR} problem with provable optimality guarantees. This paper provides an efficient algorithm to compute solutions of {\sc UMDR} problem that are close to the optimal solution, with a precise theoretical guarantee on their {\it approximation ratio} (that is, the deviation of the computed solution from the optimal solution). Moreover, the algorithm yields a scalable running time, which is necessary in event-based DR for MGs in practice.

A number of prior studies \cite{6477803, 6908041, 4494595, 6716076, 6980137, 6840288} formulated the DR management problem as an optimal power flow (OPF) problem with the objective of maximizing the cumulative user utilities or minimizing the power supply costs, losses or curtailed load, by considering continuous decision variables. If the (discrete) binary decision variables ($x_k$) in {\sc UMDR} problem are allowed to be {\it continuous} variables, such that $x_k \in [0, 1]$, then the resulting problem becomes a {\it convex quadratic programming} problem, which can be solved optimally in polynomial time, for e.g., by applying Interior Point methods. Furthermore, a slightly modified variant of the proposed algorithm with a worst-case guarantee notably close to $1$ can be devised while requiring only an order of $\mathcal{O}(n\log n)$ running time. However, setting the decision variable ($x_k$) to be discrete or continuous alters not only the complexity of {\sc UMDR} problem but also affects its practical application aspects. Particularly, the continuous case corresponds to customers having elastic (i.e., divisible) power demands allowing partial curtailment. Whereas, the discrete case requires a customer load either to be curtailed or retained holistically. This paper focuses on the latter case where binary decision variables ($x_k$) are incorporated since it is much harder to solve and is more realistic in practice.

The non-linear power and voltage constraints present in the OPF problem make it practically inefficient to exactly solve the DR management problem in distribution systems even with hundreds of customers. Most of prior work consider the convex relaxation of the problem instead or employ conventional numerical solvers to tackle the non-linear optimization problem. The conventional numerical solvers, as the one adopted in this study, usually are not guaranteed to terminate within a reasonable time frame (unless otherwise explicitly restricted) and necessarily return optimal solutions for integer programming problems. For example, we observe that Gurobi optimizer sometimes failed to terminate and return a solution to {\sc UMDR} problem considering the specified optimality parameters.


For the purpose of this study, we focus on analyzing and solving the event-based DR management problem in islanded MGs focusing on their special requirements, philosophy of operation and problem specific practical constraints. Typically, for large-scale power systems, power flow equations are included to guarantee that the voltage levels are within the IEEE standard limits. This is unlikely to be an issue for isolated MGs where the size of the grid is typically small and thus power flow equations have been excluded from {\sc UMDR} problem. It is assumed that by maintaining a balance between the generation and demand on the island the system voltages are within the acceptable range. This assumption is also followed in \cite{nikos4454000, 6665102, 6510507, 6891360}.

In order to measure the quality of the computed solution, a measure denoted by approximation ratio is defined as follows.




\medskip

\begin{definition}
Let $ x^\ast_k $ be an optimal solution to {\sc UMDR} problem and $\Opt \triangleq \sum_{k\in{\cal N}}u_k x^\ast_k$ be the corresponding total utility.
An approximate solution with  worst-case guarantee $\alpha \in [0, 1]$ to {\sc UMDR} problem is a feasible solution $(\hat{x}_k)_{k \in {\cal N}} \in \{0, 1\}^n$ satisfying: 
\begin{eqnarray}
& &   \displaystyle \sum_{k\in{\cal N}}u_k \hat{x}_k\ge \alpha \cdot \Opt \\
\text{and} & & \displaystyle \Big|\sum_{k\in{\cal N}}S_k \hat{x}_k\Big| \le C \label{C1'}.
\end{eqnarray}
\end{definition}

\medskip

The worst-case guarantee $\alpha$ is also called the {approximation ratio}, which characterizes the ratio between the optimal solution and the approximate solution. When $\alpha = 1$, this becomes an exact optimal solution. In the subsequent sections, an efficient algorithm is presented with a definite worst-case guarantee $\alpha$ independent of the number of customers.



\section{Efficient Algorithms}\label{sec:algs}

In this section, three greedy algorithms are presented, where the first two follow from the conventional strategies, and the last one is a novel strategy proposed in this paper. Without loss of generality, assume $|S_k| \le C$ for all $k$ (since any demand $k$ with $|S_k|>C$ cannot be part of any feasible solution). An allocation $(x_k)_{k \in {\cal N}}$ can be equivalently represented by the set of satisfied customers $X \triangleq \{k\in {\cal N}\mid  x_k=1\}$. For a subset $X \subseteq {\cal N}$, denote $u(X) \triangleq \sum_{k\in X}u_k$. 

\begin{figure}[!htb]
\centering
 \includegraphics[scale=0.6]{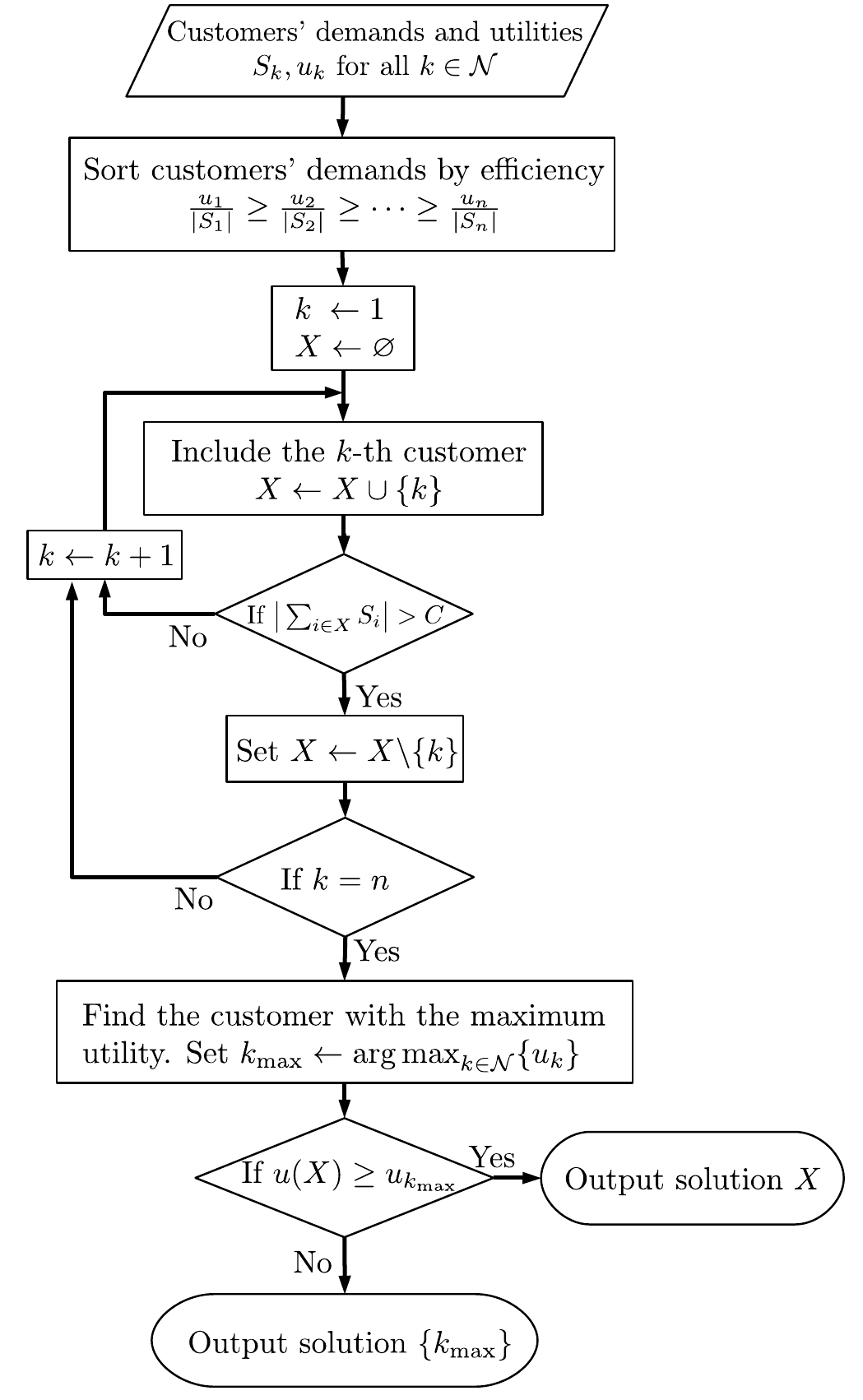}
\caption{Flow chart for Greedy Ratio Algorithm {\sc (GRA)}.} \vspace{-15pt} 
\label{fig:alge}
\end{figure}

The three greedy algorithms are described as follows.
\begin{enumerate}

\item
{\em Greedy Utility Algorithm} {\sc(GUA)}: First, sort the customers in $\cN = \{1,...,n\}$ by their utilities in a non-increasing order (with arbitrary tie-breaking), such that 
\begin{equation}
u_1 \ge u_2 \ge ... \ge u_n \,.
\end{equation}
Then, select the satisfiable demands sequentially from the first customer according to the order whenever feasible (i.e. $\big|\sum_{k}S_k x_k\big| \le C$).

\item 
{\em Greedy Demand Algorithm} {\sc(GDA)}: Similar to {\sc GUA}, but sort the customers by the magnitudes of their demands in a non-decreasing order, such that  
\begin{equation}
|S_{1}| \le |S_{2}| \le ... \le |S_{n}| \,.
\end{equation}
Then, select the satisfiable demands sequentially from the first customer according to the order whenever feasible (i.e. $\big|\sum_{k}S_k x_k\big| \le C$).

\item 
{\em Greedy Ratio Algorithm} {\sc (GRA)}: Similar to {\sc GUA} and {\sc GDA}, but sort the customers by the efficiency (i.e., $\frac{u_k}{|S_k|}$) in a non-increasing order, such that  
\begin{equation}
\label{eqn:gra} 
\frac{u_1}{|S_1|} \ge \frac{u_2}{|S_2|} \ge ... \ge \frac{u_n}{|S_n|} \,.
\end{equation}
Then, obtain the greedy solution $X$ by selecting sequentially from the first customer according to the order whenever feasible (i.e. $\big|\sum_{k\in X}S_k x_k\big| \le C$). Find the single customer with the highest utility ($\max_{k\in \cN}\{ u_k\}$). Lastly, output the solution with the higher utility between $X$ and the highest utility customer. A flowchart for {\sc GRA} is depicted in Fig.~\ref{fig:alge}.

\end{enumerate}


It is worthy to note that {\sc GUA} and {\sc GDA} are the common strategies of standard load curtailment in practice \cite{5348255}. However, unlike {\sc GRA}, they may produce arbitrarily worse solutions for {\sc UMDR} problem when compared to the optimum solution. This claim could be easily verified by considering the following two simple examples of {\sc UMDR} problem. 
	
\begin{example} Consider an instance with $n=2$ and $C=F$, where $F$ is a real positive number. For simplicity, consider setting zero reactive power for all customers. Customer $1$ has demand of $|S_1| = 1$ and utility of $u_1 = 1$ and customer $2$ has demand of $|S_2| = F$ and utility of $u_2=F$.
\end{example}
	
The optimal solution set would consist of customer $2$ thus reaching an optimal solution value of $F$. Whereas, {\sc GDA} would select customer $1$ thus producing a solution with value of $1$. For a suitably large selection of $F$ the relative performance guarantee of {\sc GDA} compared to the optimal solution could be arbitrarily close to $0$. 
\begin{example} Consider the following instance with $n=b$ and $C=10$, where $b$ is an integer positive number greater than $2$. For clarity, consider setting zero reactive power for all customers. Customer $1$ has a demand of $|S_1|=10$ and utility of $u_1=10$ and customers $2$ to $b$ have demand of $10/b$ and utility of $9$.
\end{example}

The optimal solution would be to select customers $2$ to $b$ thus producing an optimal solution with value $9(b-1)$. {\sc GUA}, instead, will select only customer $1$ thus generating a solution with a value of only $10$. For considerably large values of $b$ the ratio between solutions generated by {\sc GUA} and optimal solution value could be arbitrarily close to $0$.

In fact, among the three algorithms, only {\sc GRA} possesses a worst-case guarantee. More precisely, {\sc GRA} provides a solution that in the worst possible case is at least $\frac{1}{2} \cos\tfrac{\theta}{2}$ of the optimal solution for {\sc UMDR} problem, as described in the following theorem. Each customer's power demand has a phase angle between the reactive and apparent powers. Define $\theta$ to be the maximum difference between the phase angles of any pair of power demands of the customers.

\medskip

\begin{customthm}{1} \label{thm:alg-greedy}
Algorithm {\sc GRA} produces a feasible solution within a worst-case guarantee of $\frac{1}{2} \cos\tfrac{\theta}{2}$ of any optimal solution of \textsc{UMDR} problem.
\end{customthm}

\smallskip

The appendix provides a detailed proof of Theorem \ref{thm:alg-greedy}. 
The worst-case guarantee of {\sc GRA} depends on the angle $\theta$; the smaller  $\theta$ is, a  better approximation guarantee is achieved. If $\theta = 0$, then all customers have the same ratio $\frac{P_k}{Q_k}$. In  practical settings, the load power factor usually varies between $0.8$ to $1$ \cite{1339347} and thus the maximum phase angle difference $\theta$ between any pair of demands is restricted to be in the range of $[0, 36^{\circ}]$. This implies that the approximation  ratio is at least $0.475$. 

It is worth mentioning that {\sc GRA} achieves the aforementioned approximation ratio only in order of $\mathcal{O}(n\log{n})$ time. One can notice from Fig. \ref{fig:alge} that the algorithm requires three basic procedures to terminate, namely sorting followed by a linear iteration over the newly sorted set and determining the customer with the highest utility. Each of the latter two procedures require $\mathcal{O}(n)$ operations, while sorting could be performed in order of $\mathcal{O}(n\log{n})$ time by adopting a well-known sorting algorithm (e.g. Merge sort). Therefore, the running time complexity of {\sc GRA} is $\mathcal{O}(n\log{n} + 2n) = \mathcal{O}(n\log{n})$.

\section{System and Settings Under Study}\label{sec:sims}

To evaluate the performance of the proposed algorithm, an MG with an overall capacity of $2$MVA and over thousands of customers is considered. Each customer has a specific power demand (including both active and reactive power) and a utility that is generated according to a probability preference model. The amount of generation on an MG is typically less than the amount of demand and thus, the customers may suffer from a reduction of generation capacity occasionally. Various types of loads are considered including residential and industrial customers ranging between $300$KVA to $1$MVA. It is also assumed that the MG is equipped with a two-way communication infrastructure capable of sending the optimal load management signals (determined by the centralized controller) and allowing for customer demand and utility to be sent back to the centralized controller. The central controller is assumed to have full control over the on/off operations of its customers.

In order to quantify the performance of the proposed algorithm, Gurobi optimizer is employed to obtain the close-to-optimal solutions numerically. In this paper, the output solution of Gurobi optimizer is denoted by $\Opt$. The following parameters were set in Gurobi optimizer: (1) the total time expended for solving the problem was $200$ seconds, (2) absolute mixed integer programming (MIP) optimality gap (i.e., the threshold of the absolute gap between the lower and upper objective bound) was set to zero, and (3) infeasibility tolerance was set to $10^{-9}$. It is worth noting that there are no guarantees that given an integer programming problem the optimizer will return an optimal solution nor it will terminate in a reasonable time (i.e., within $200$ seconds for each run). Whenever the optimizer exceeds the time limit, the current best solution is considered to be optimal.

The simulations were evaluated using 2 Quad core Intel Xeon CPU E5607 2.27 GHz processors with 12 GB of RAM. The algorithms were implemented using Python programming language with Scipy library for scientific computation. The power system analysis package PSCAD was utilized for testing the load curtailment schemes produced by the proposed algorithm. Typically, the load power factor varies between $0.8$ to $1$ and thus the maximum phase angle $\theta$ between any pair of demands is restricted to be in the range of $[0, 36^{\circ}]$.

Various case studies are performed to evaluate the proposed algorithm by taking into account the correlation between customer load and utility considering various load types. The following are settings for the case studies in this paper.
 
\begin{enumerate}

\item[(i)] {\em Utility-demand correlation}:
\begin{enumerate}

\item {\em Correlated setting (C)}: The utility of each customer is a function of the power demand:
\begin{equation}\label{eq:valuationfunction}
u_k({|S_k|}) = a\cdot {|S_k|}^2 + b\cdot {|S_k|} + c \,,
\end{equation}
where $a > 0, b, c \ge 0$ are constants. For simplicity, $u_k(|S_k|) = |S_k|^2$ is considered in the simulation.

\item {\em Uncorrelated setting (U)}: The utility of each customer is independent of the power demand and is generated randomly from $[0, |S_{\max}(k)|]$. Here $S_{\max}(k)$ depends on the customer type (as defined below): if customer $k$ is an industrial customer then $|S_{\max}(k)| = 1$MVA, otherwise $|S_{\max}(k)| = 5$KVA.

\end{enumerate}

\item[(ii)] {\em Customer types}:
\begin{enumerate}

\item {\em Residential (R) customers}: The customers are comprised of residential customers having small power demands ranging from 500VA to 5KVA.

\item {\em Mixed (M) customers}: The customers are comprised of a mix of industrial and residential customers. Industrial customers have big power demands ranging from 300KVA up to 1MVA and constitute no more than 20\% of all customers chosen at random.

\end{enumerate}

\end{enumerate}

In this paper, the case studies will be represented by the aforementioned acronyms. For example, the case study named CM stands for the one with mixed customers and utilities-demand correlation.

%
%
%
%
%
%

\section{Simulation Results}
\label{simresutssec}
\subsection{Comparative Analysis}

In this subsection the proposed approaches are compared in terms of quality of solution. The optimal solutions computed by Gurobi are considered to be the base case for the comparison. The algorithms are applied to various case studies where each case study is analyzed considering changes in the set of customers. As an example, {\sc GRA} is applied 30 times for each of the $m$ number of customers (where $m$ varies between $100$ to $1500$ in steps of hundred) for case study CR (i.e., correlated, residential) considering random changes in demands and utilities of customers. Thus, the total number of experiments for each case study is $450$. In particular,  Table~\ref{tab:1} highlights the results obtained using {\sc GRA}, {\sc GUA}, and {\sc GDA} for the various case studies. The results in Table~\ref{tab:1} present the minimal ratio between the solutions obtained by the proposed algorithms and Gurobi. It is worthy to note that the closer this value is to $1$, the closer is the solution to the optimum.



\begin{table}[h]
	\caption{The worst-case approximation ratios of {\sc GRA}, {\sc GUA}, and {\sc GDA} with respect to the optimal solutions computed by Gurobi optimizer.}
	\label{tab:1}
	\renewcommand{\arraystretch}{1.4}
	\centering
	\begin{tabular}{cc|c|c|}
		\cline{3-4}
		& \textbf{} & \textbf{C} & \textbf{U}  \\ \hline
		\multicolumn{1}{|c|}{\multirow{2}{*}{\textbf{GRA}}} & \textbf{R} & 0.999 &  0.883  \\ \cline{2-4} 
		\multicolumn{1}{|c|}{} & \textbf{M} & 0.921 & 0.568  \\ \hline
		\multicolumn{1}{|c|}{\multirow{2}{*}{\textbf{GUA}}} & \textbf{R} & 0.999 & 0.934  \\ \cline{2-4} 
		\multicolumn{1}{|c|}{} & \textbf{M} & 0.921 & 0.403  \\ \hline
		\multicolumn{1}{|c|}{\multirow{2}{*}{\textbf{GDA}}} & \textbf{R} & 0.541 & 0.839  \\ \cline{2-4} 
		\multicolumn{1}{|c|}{} & \textbf{M} & 0.002 & 0.262   \\ 
		 \hline
	\end{tabular}

\end{table}

\begin{figure*}
	\centering
	\begin{minipage}{.49\textwidth}
		\centering
		\includegraphics[scale=0.45]{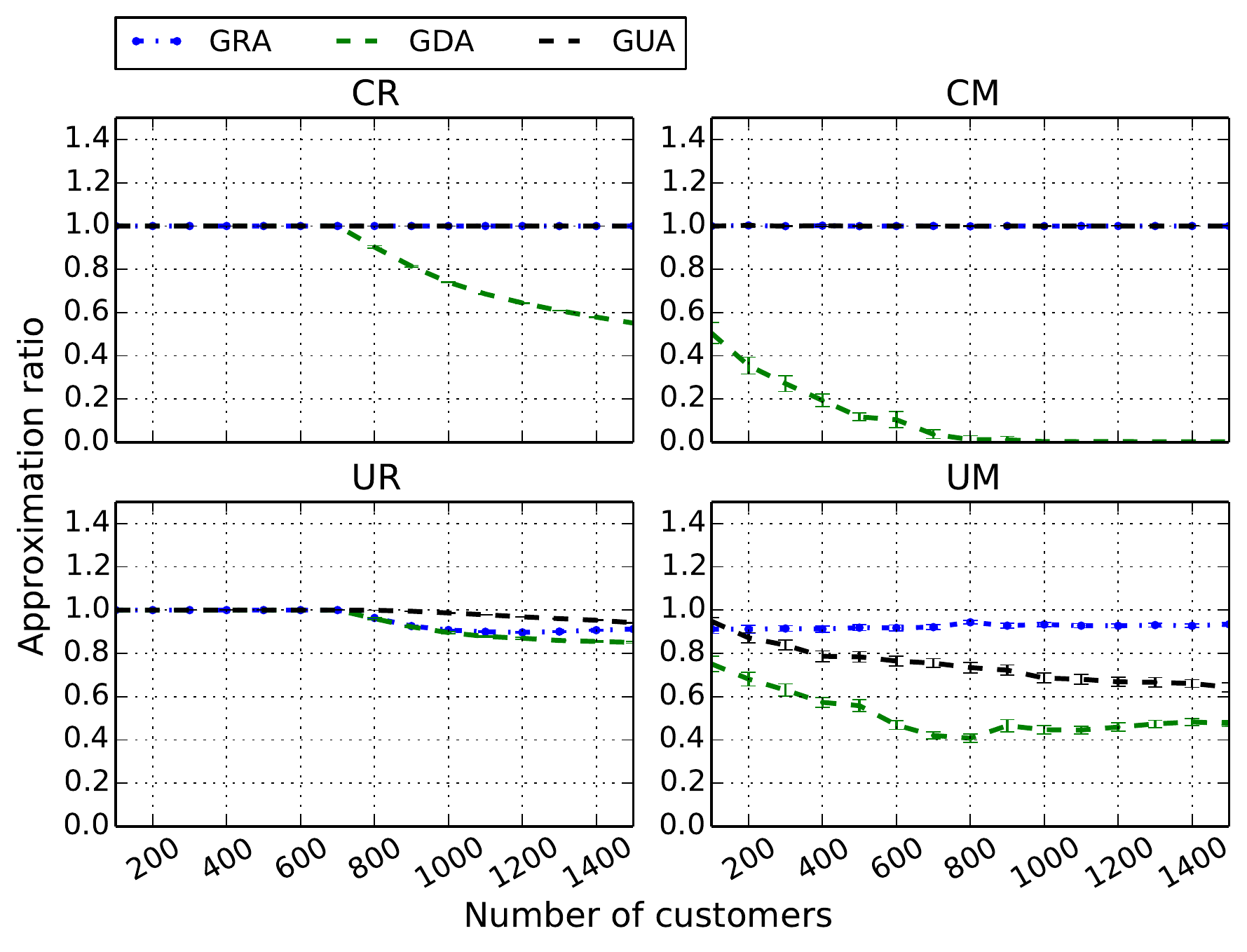}
		\captionof{figure}{The average approximation ratios of  {\sc GRA}, {\sc GUA} and {\sc GDA} against the number of customers at 95\% confidence interval.}
		\label{fig:cr}
	\end{minipage}%
	\hspace*{5pt}
	\begin{minipage}{.49\textwidth}
		\centering
		\includegraphics[scale=0.45]{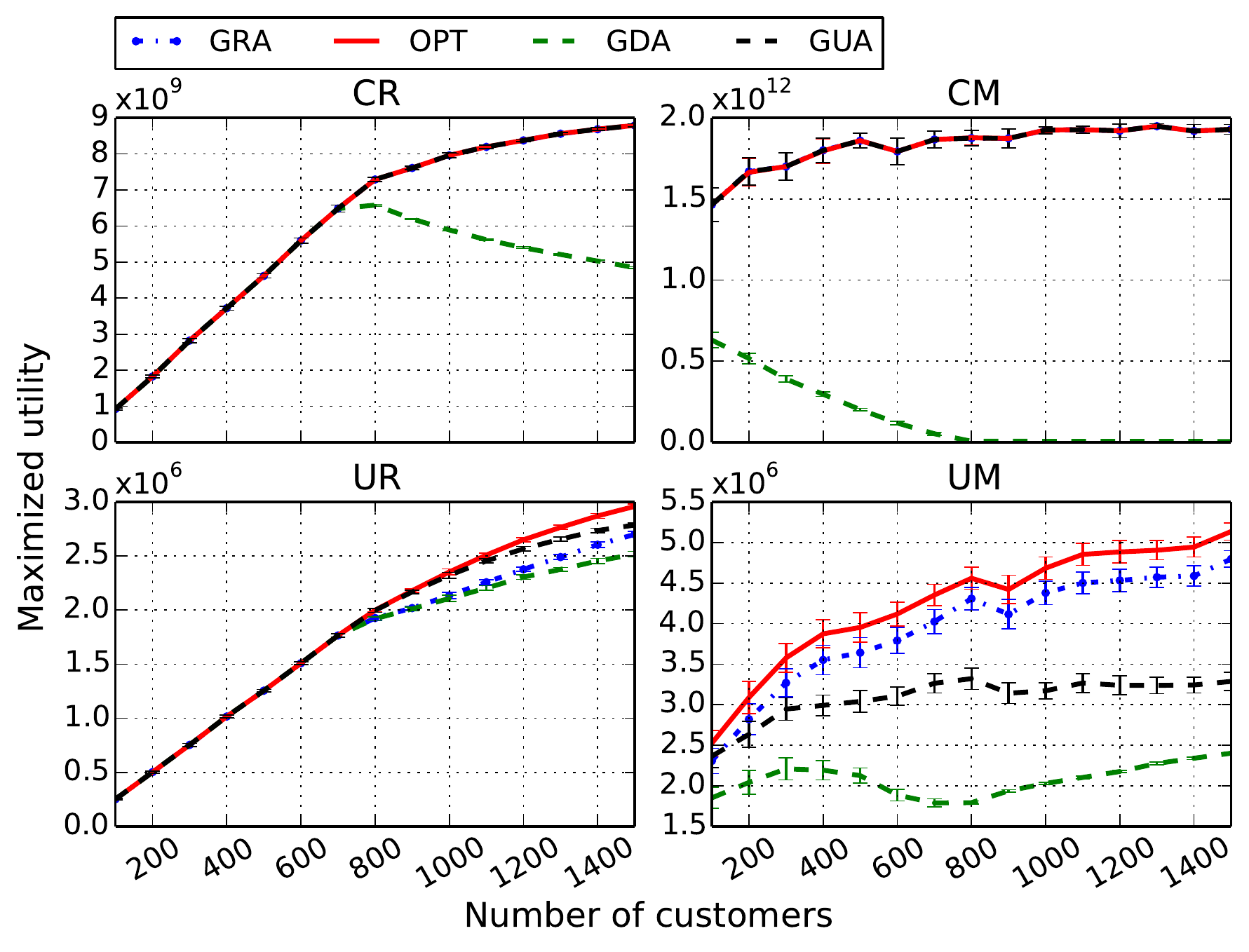}
		\captionof{figure}{The average objective values of {\sc GRA}, {\sc GUA}, {\sc GDA} and $\Opt$  against the number of customers at 95\% confidence interval.}
		\label{fig:obj}
	\end{minipage}
	\vspace*{-15pt}
\end{figure*}

As can be observed from Table \ref{tab:1} when there is a utility-demand correlation for all customers, {\sc GRA} and {\sc GUA} obtain solutions that are very close to the optimal $\Opt$. Furthermore, for all the case studies {\sc GRA} provides the highest approximation ratio in the worst-case (worst-case of {\sc GRA} is $0.568$, while that of {\sc GUA} is $0.403$) when compared to other candidate algorithms. Figure~\ref{fig:obj} presents the maximized utility for the three algorithms at $95$\% confidence interval considering different customer set cardinality. Similarly, each point presented in Fig.~\ref{fig:obj} represents the average objective value among $30$ iterations. It was observed experimentally that $30$ iterations was sufficient for convergence of the sample mean and variance.

It is observed that {\sc GDA} performs the worst in terms of the quality of solution. In fact, {\sc GDA} performance degrades significantly as the number of customers increase when there is correlation between utilities and demands. The reason is that {\sc GDA} heuristically selects the smallest demands first which on contrary obtain lower utility when considering a quadratic utility function. This situation does not occur when there is no such correlation (namely, in case studies UR and UM). For case study UR, the observed performance of {\sc GUA} is slightly better than that of {\sc GRA} when considering large customer participation. This could be due to the fact that with increasing customer participation the probability of having a customer with a small utility but high utility to demand ratio increases. Unlike {\sc GUA}, {\sc GRA} selects those customers first thus failing to consider the customers with relatively large utilities. All algorithms can output optimal solutions when having only few residential customers, because at this scale all customers' demands are below the total capacity of $2$MVA and hence can be all satisfied (see case studies CR and UR in Fig.~\ref{fig:obj}).

The observed approximation ratios are plotted in Fig.~\ref{fig:cr} against the number of customers for each scenario. When a curve is close to the line $y=1$, it indicates a close to optimal solution. The observed approximation ratios of {\sc GRA} and {\sc GUA} are approximately 1 for case studies CR and CM for scalable number of customers. 
As highlighted earlier, the theoretical guarantees on the worst-case approximation ratio of {\sc GRA} is $0.475$, computed using $\theta=36$ degrees (refer to Theorem \ref{thm:alg-greedy}). Nevertheless as can be observed from Table ~\ref{tab:1} {\sc GRA}, for the majority of cases, can achieve almost twice higher approximation ratio compared to the theoretical bound.

\vspace*{-6pt}
\subsection{Computational Running Time}

As stated earlier a major problem with most of the available centralized demand management schemes is the computational time. It is very important to develop fast and efficient algorithms capable of executing optimal decisions when considering significant customer participation. Thus, one of the key parameters utilized to evaluate the performance of the proposed algorithm is the computational time. 

In this subsection the computational time of the proposed {\sc GRA} is compared against the Gurobi solver. Computational time is of significant importance when designing centralized controllers for MGs since this will have implications on the stability of MG. Note that the running time complexity of {\sc GUA} and {\sc GDA} is the same as that of {\sc GRA}. For clarity of presentation, however, we investigate only the running time of {\sc GRA}.

\begin{figure}[!htb]
	\begin{center}
		\includegraphics[scale=.45]{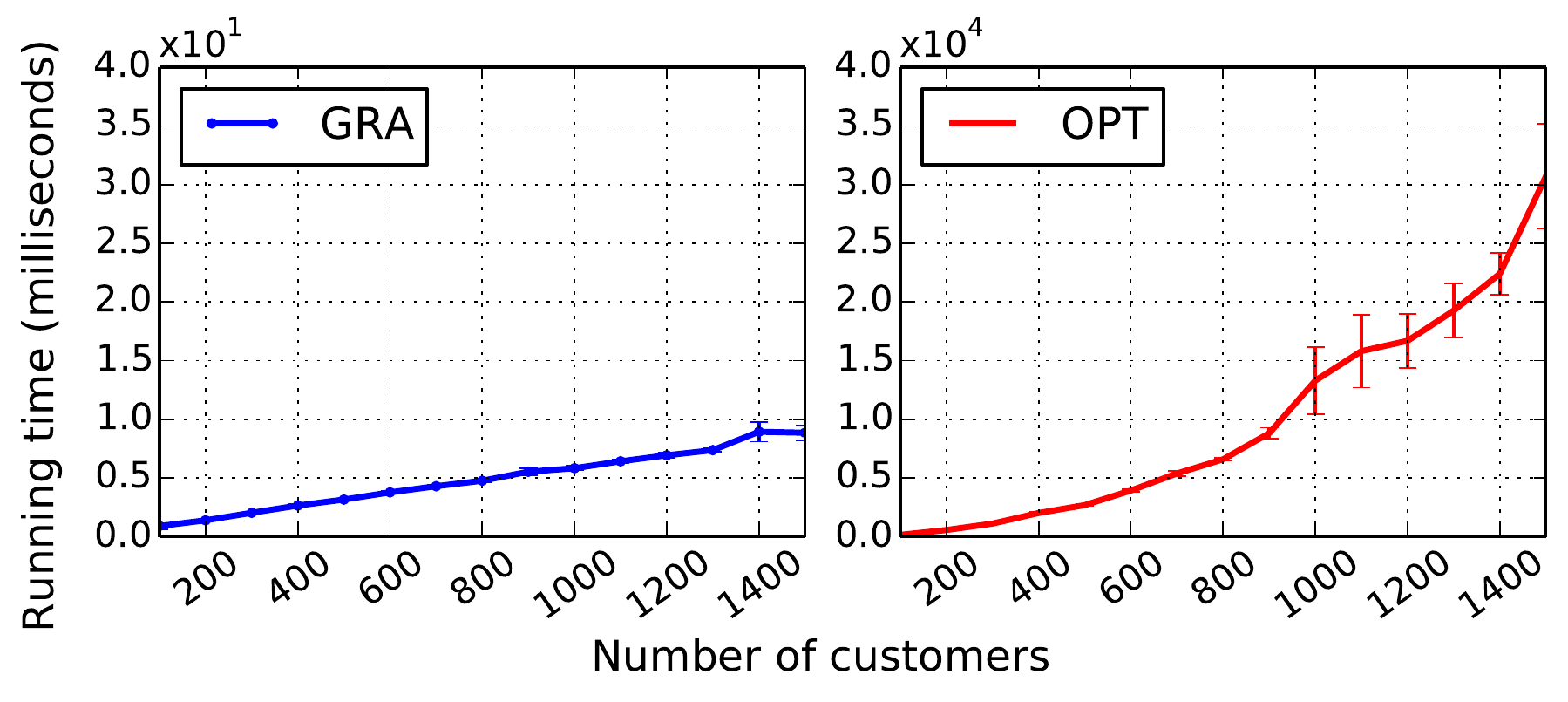}
	\end{center}
	\caption{The average running time of {\sc GRA} (left) and $\Opt$ (right) against the number of customers at 95\% confidence interval.}
	\label{fig:time}
\end{figure}

The running time of {\sc GRA} is compared to that of Gurobi optimizer in Fig.~\ref{fig:time}. For brevity, only the results related to a case study CR are presented in Fig.~\ref{fig:time}. However, it is worthy to note that nearly the same running time was observed for {\sc GRA} and Gurobi optimizer when considering the rest case studies. The computational time reported is the average running time over $30$ iterations. As can be observed, for an MG with roughly $600$ customers the centralized controller operated using Gurobi solver will take roughly $5$ seconds. 

In these case studies it is expected that the MG will not be capable of maintaining stable operation. On the contrary, for the same number of customers the computational time needed for the {\sc GRA} is nearly $5$ milliseconds. This would allow load management decisions to be made almost instantaneously and thus allowing for the MG to stabilize. 

For an MG with up to $1400$ customers, the {\sc GRA}  is capable of providing close to optimal solutions in less than $10$ milliseconds. Furthermore, it is worthy to mention that for some cases the Gurobi optimizer did not return an optimal solution within $200$ seconds, but on the contrary {\sc GRA} can always output a solution in nearly linear time. 

%
\vspace*{-25pt}
\subsection{Dynamic Generation Capacity}

The previous subsections considered an MG with a fixed generation capacity of $2$MVA and varying load capacities. Here simulations are performed considering the case when the MG's generation capacity is varying over time due to possible occurrences (e.g., due to failure, maintenance or resumption), whereas the set of customers of cardinality $1000$ is fixed. 

\begin{figure}[!htb]
	\begin{center}
		\includegraphics[scale=0.45]{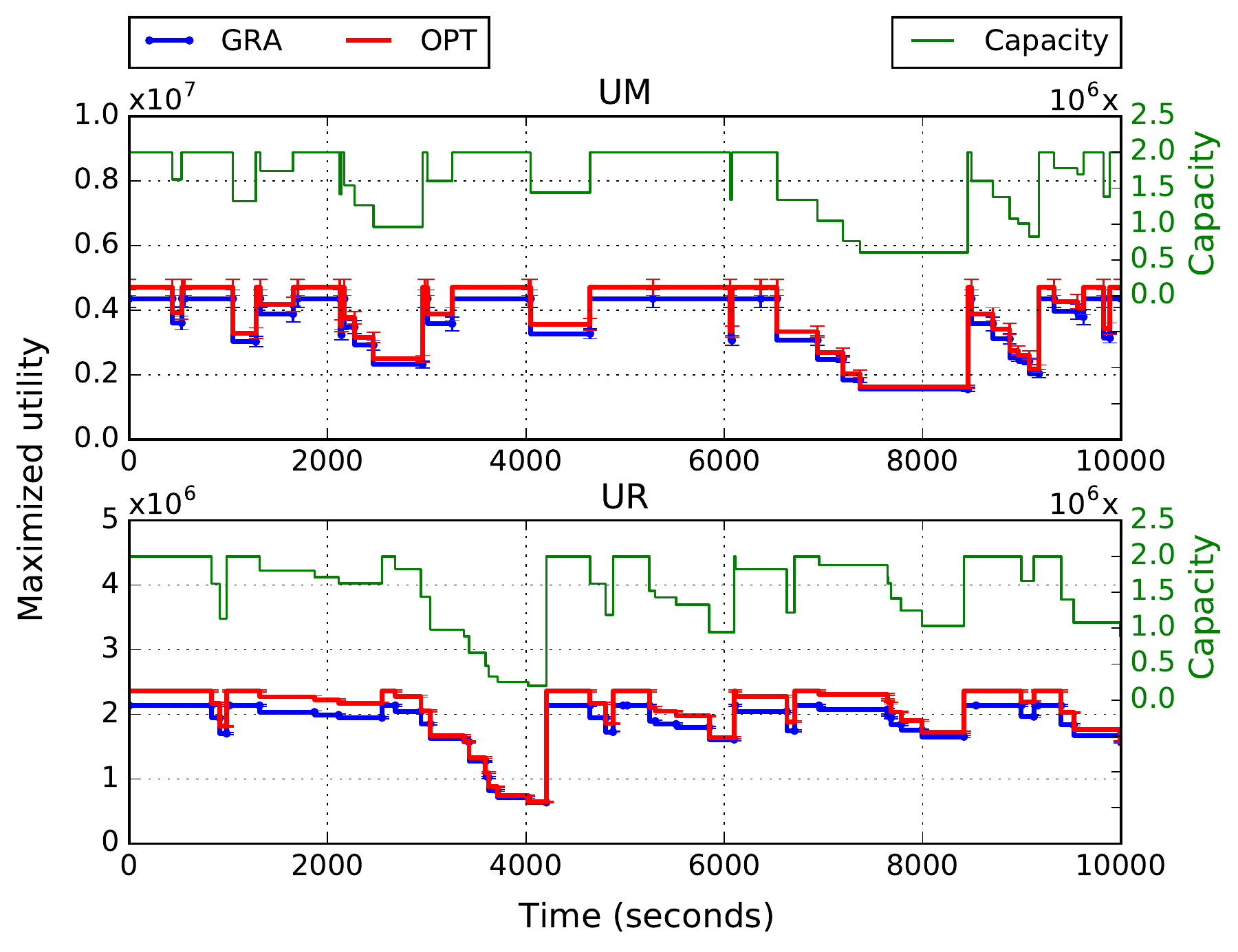}
	\end{center}
	\caption{The maximum utility of {\sc GRA} and $\Opt$ considering dynamic generation capacity.}
	\label{fig:dyn}
\end{figure}

The generation capacity of MG is dynamically varied between $100$KVA and $2$MVA and is applied to the proposed problem formulation from time $0$ to $10000$ (seconds). The events, namely Failure and Resumption, occur according to an exponential distribution with a rate of $0.005$. When the MG is in the Failure state, the generation capacity decreases randomly from $5\%-35\%$, whereas when in the Resumption state the MG's generation capacity is fully resumed (i.e.,  $C=2$MVA). Whether an event is in a Failure or Resumption state, is determined according to a Markov chain with the following settings:

\begin{enumerate}
	\item Steady $\rightarrow$ Failure with a probability of $65\%$.
	\item Steady $\rightarrow$ Resumption with a probability of $35\%$.
\end{enumerate}

 
In Fig.~\ref{fig:dyn}, {\sc GRA}'s maximum utility is slightly lower than $\Opt$ in both case studies UM and UR (other case studies, namely CM and CR, {\sc GRA} is even closer to $\Opt$). It is observable that the gap between {\sc GRA} and $\Opt$ is reduced when the generation capacity is reduced due to the Failure events. The results show that {\sc GRA}  is capable of obtaining close to optimal solutions under various case studies with minimal computational time. 

 

 
\subsection{Application on a Feeder of Benchmark Canadian System}
\label{sec:appfeeder}
\begin{figure}[!htb]
	\begin{center}
		\includegraphics[scale=.43]{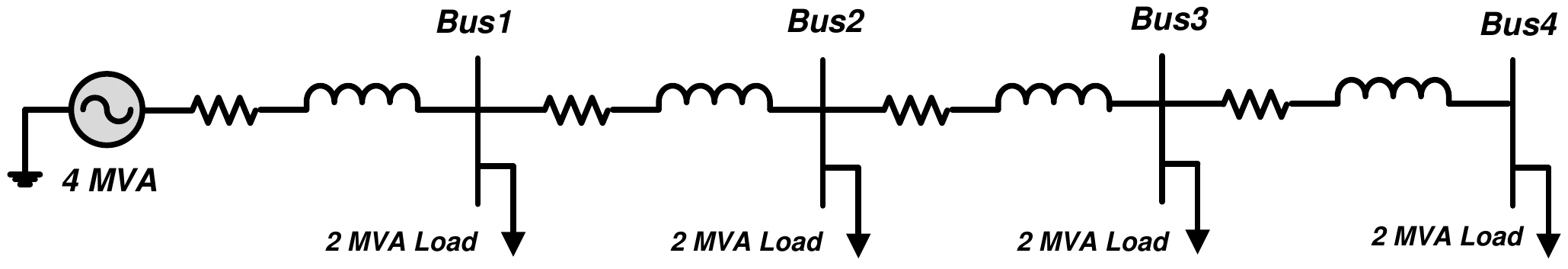}
	\end{center}
	\caption{A 4-bus feeder from Canadian benchmark distribution system.}
	\label{fig:system}
\end{figure}

The proposed algorithm is applied to one of the feeders of the Canadian benchmark bystem. The power system simulation software PSCAD has been employed to simulate the feeder which is rated at $8.7$MVA, $400$A and $12.47$KV. Each feeder section is a $700$MCM Cu XLPE cable with $Z=0.1529 + J 0.1406$ $\Omega/km$. Each node consists of a $2$MVA total load as given in Fig. \ref{fig:system}. Overall, up to $1500$ customers have been allocated among the nodes. 

\begin{figure*}
	\centering
	\begin{minipage}{.49\textwidth}
		\centering
		\includegraphics[scale=.45]{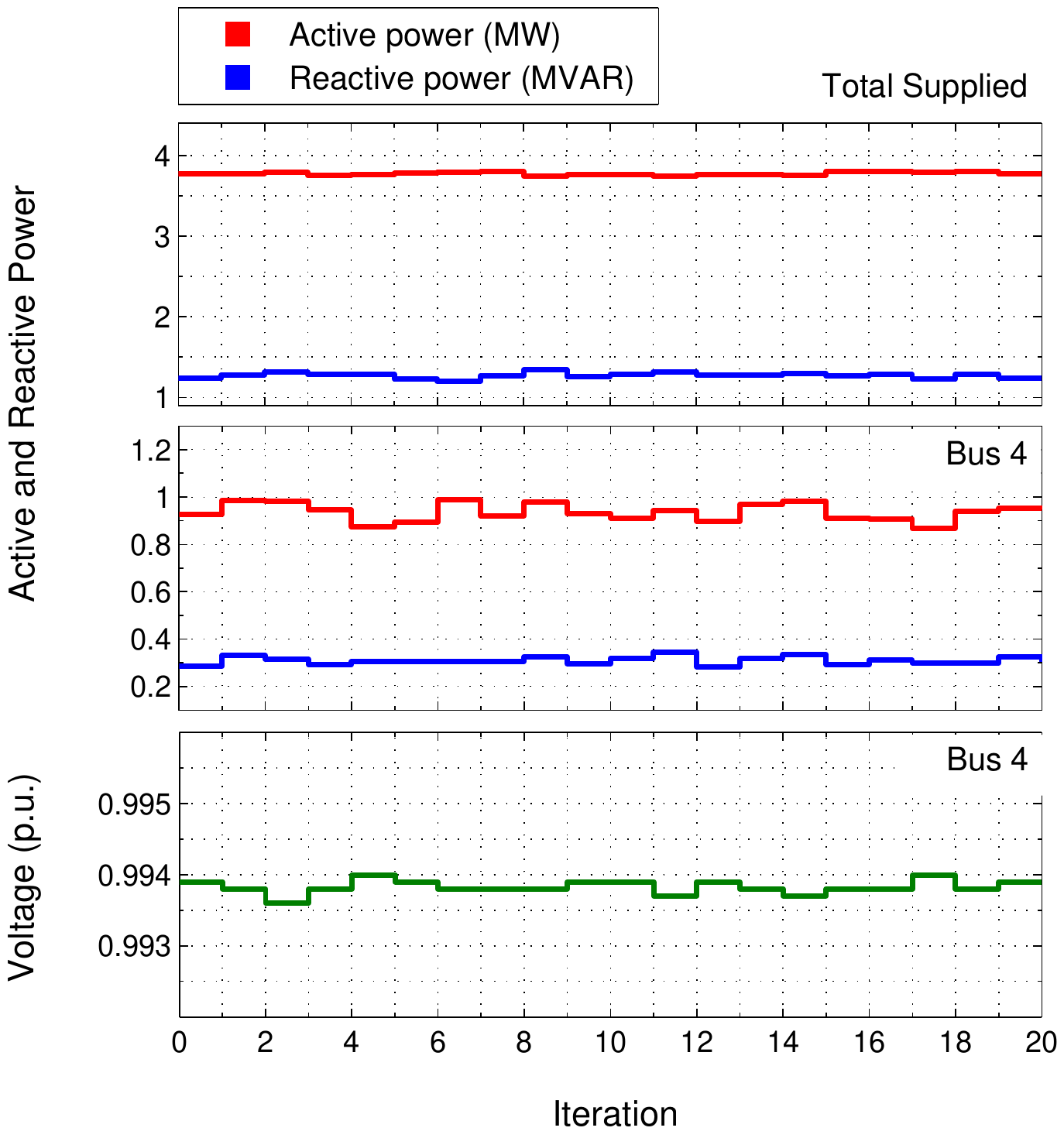}
		\captionof{figure}{The generation capacity of MG, load and voltage profile on bus 4 for case study CR.}
		\label{fig:fcrw}
	\end{minipage}%
	\hspace*{5pt}
	\begin{minipage}{.49\textwidth}
		\centering
		\includegraphics[scale=.45]{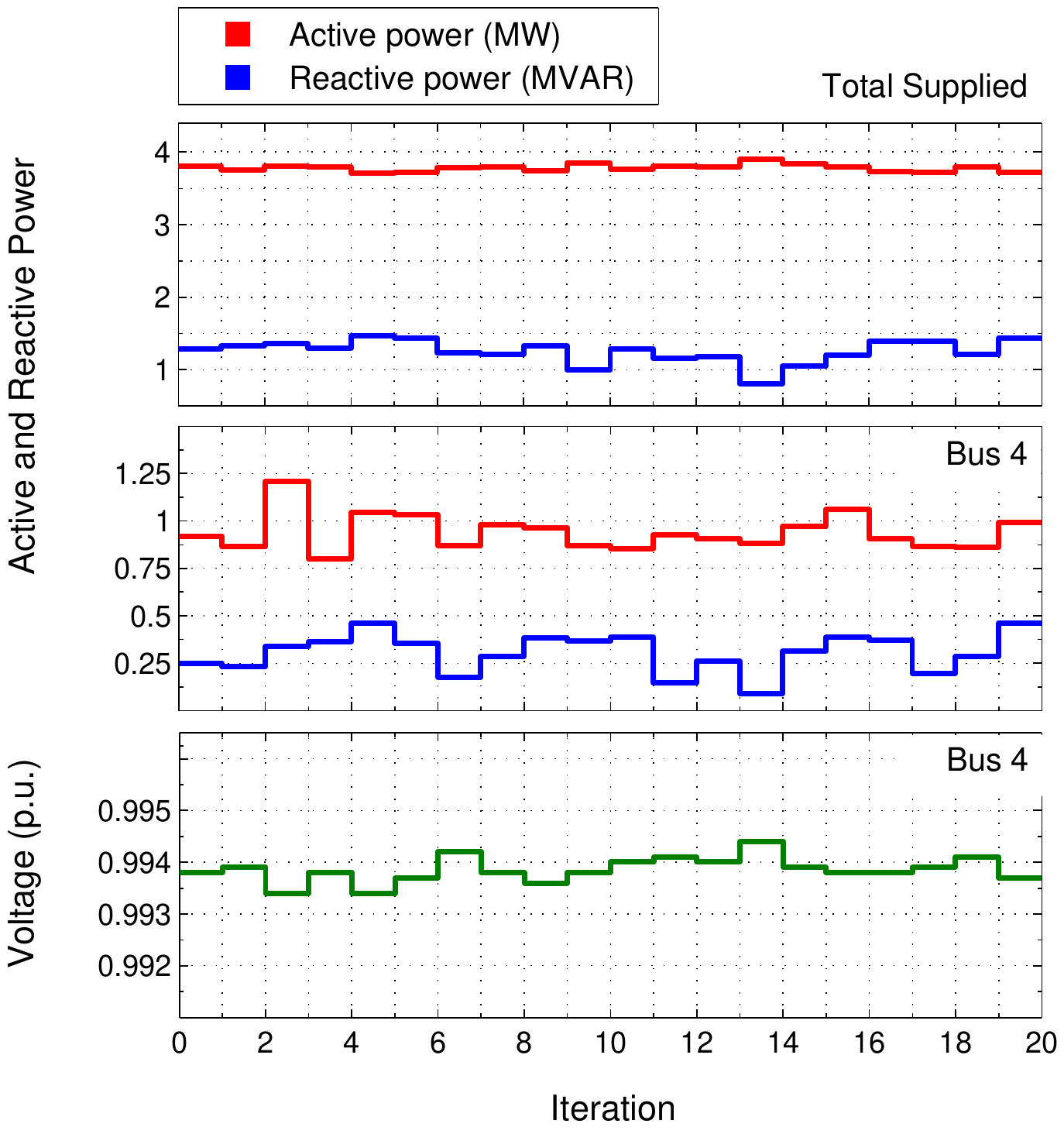}
		\captionof{figure}{The generation capacity of MG, load and voltage profile on bus 4 for case study CM.}
		\label{fig:fcmw}
	\end{minipage}
	\vspace*{-15pt}
\end{figure*}

To account for the losses in the system, the worst case is executed where it is assumed that a total of $8$MVA will be allocated at the end of feeder (bus 4). The maximum losses expected under such case is  $P_{loss}=24.4$KW and $Q_{loss} = 22.496$KVA. This amount of power losses is allocated in each iteration and withheld from the maximum MG DG capacity. This is done to assure sufficient power balance between load, losses and generation. Figures \ref{fig:fcrw} and \ref{fig:fcmw}  present the active and reactive demand generated by the DG source, consumed by the load at bus 4 as well the voltage for 20 iterations of utility bidding for case studies CR and CM respectively. As can be seen, the total generated capacity is less than $4$MVA to account for losses. The voltage of the end of the feeder (bus 4) is within the IEEE standard 1547 limits.

%

	\vspace*{-25pt}
	\subsection{Network and Voltage Constraints}
	\label{graopfalg}
\begin{figure}[!htb]
			\begin{center}
				\includegraphics[scale=.45]{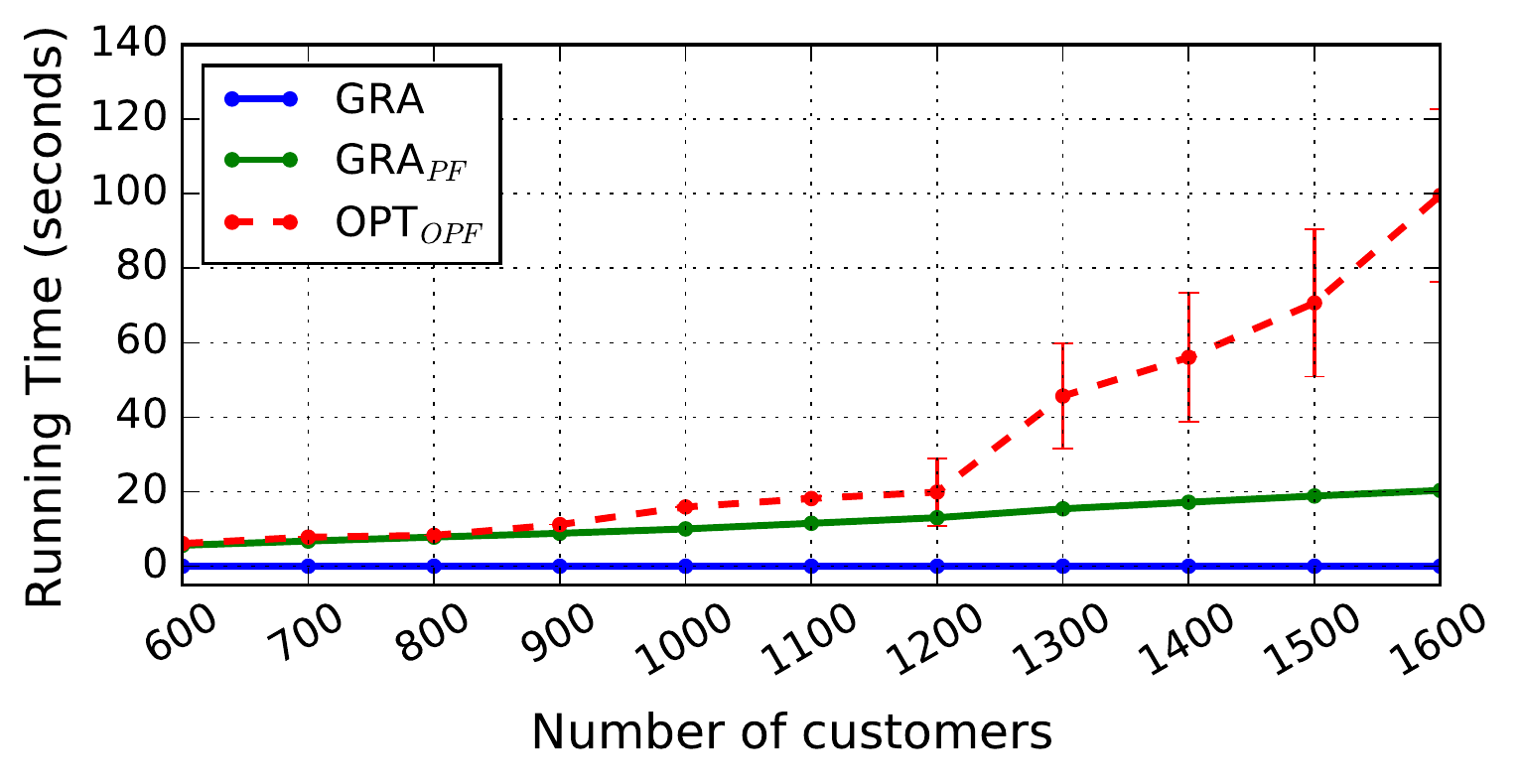}
			\end{center}
			\caption{The average running time of {\sc GRA$_{\textsc{PF}}$}, {\sc GRA} and {\sc OPT$_{\textsc{OPF}}$}, against the number of customers at $95$\% confidence interval.}
			\label{fig:nnt}
		\end{figure}

	In this subsection the proposed algorithm {\sc GRA} is extended to consider network and voltage constraints, and is applied on a feeder from the Canadian benchmark distribution system depicted in Fig.~\ref{fig:system}. The extended algorithm, denoted by {\sc GRA$_{\textsc{PF}}$} is a heuristic approach for producing load curtailment schemes which maximize the cumulative utility of the satisfied customers without violating the network and voltage constrains. The execution of {\sc GRA$_{\textsc{PF}}$} is identical to that of {\sc GRA} but with an additional subroutine of performing load flow and checking the capacity, voltage and network constrains with each customer addition. However, the feasibility of the generated solutions by {\sc GRA$_{\textsc{PF}}$} comes with a cost of increased running time (refer to Fig.~\ref{fig:nnt}) and the absence of theoretical guarantees as compared to {\sc GRA}. 
	\vspace*{-10pt}

	Furthermore, the {\sc UMDR} problem is formulated in an optimal power plow framework, and is denoted by {\sc UMDR-OPF} which is presented in the appendix. To compare the quality of solution and computational running time of both {\sc GRA} and {\sc GRA$_{\textsc{PF}}$} algorithms, we generalize the previous simulations in subsection~\ref{sec:appfeeder}. The solutions computed by Gurobi optimizer for {\sc UMDR-OPF} problem denoted by {\sc OPT$_{\textsc{OPF}}$} are considered as a base case for comparison. As mentioned in subsection~\ref{sec:appfeeder}, the amount of $P_{loss}=24.4$KW and $Q_{loss} = 22.496$KVA is suppressed from MG generation capacity in each iteration when {\sc GRA} is executed. On the other hand, for {\sc GRA$_{\textsc{PF}}$}, the problem is solved considering full MG generation capacity (i.e., $C=4$MVA) since the load flow equations account for power losses.
		\begin{figure}[!htb]
			\begin{center}
				\includegraphics[scale=.45]{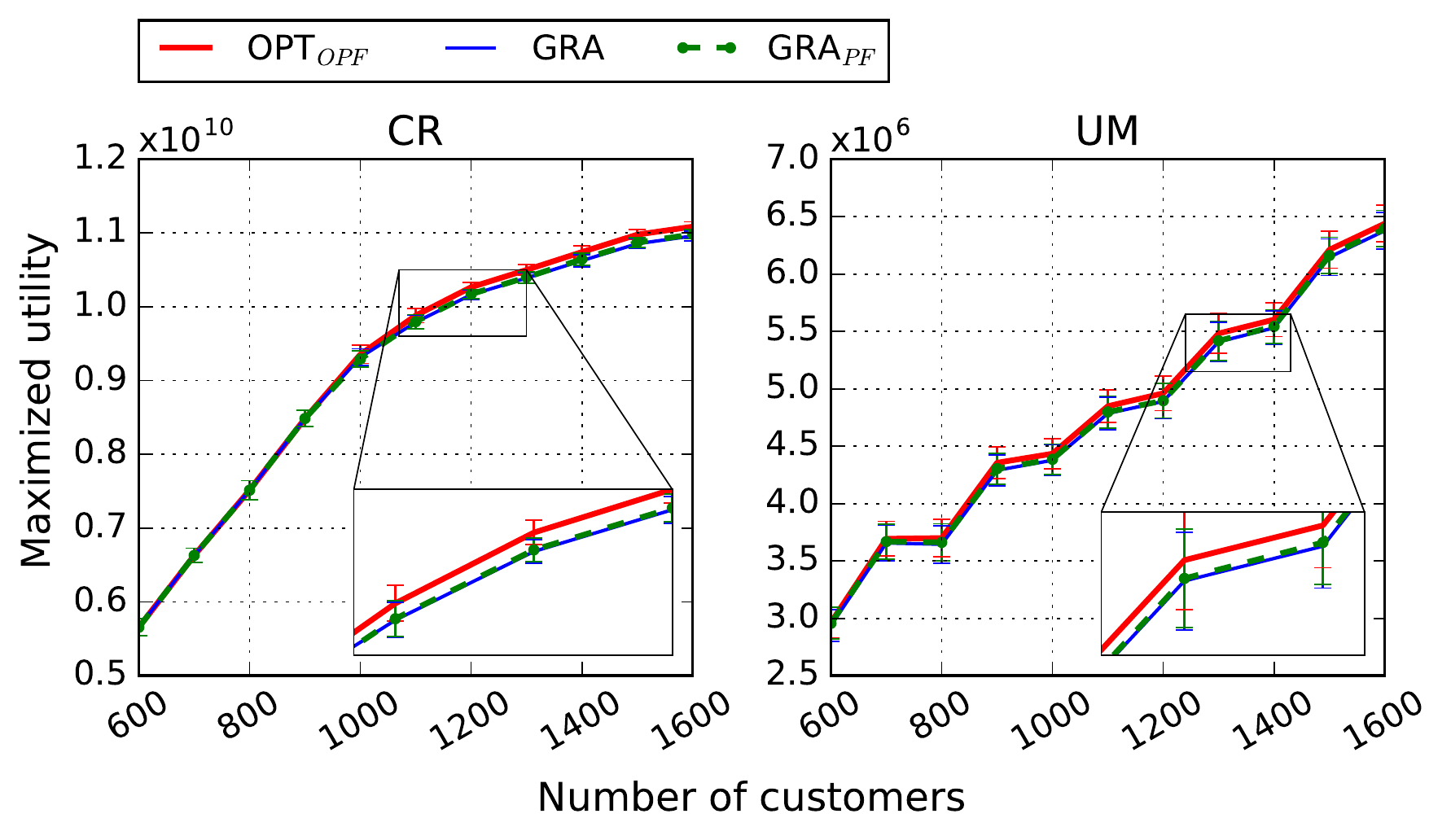}
			\end{center}
			\caption{The average objective values of {\sc GRA$_{\textsc{PF}}$}, {\sc GRA} and {\sc OPT$_{\textsc{OPF}}$} against the number of customers at $95$\% confidence interval.}
			\label{fig:nn}
		\end{figure}
		\vspace*{-15pt}
	
	 The results in Figs~\ref{fig:nn} and \ref{fig:nnt} illustrate the maximized utility and running time of {\sc GRA}, {\sc GRA$_{\textsc{PF}}$} and {\sc OPT$_{\textsc{OPF}}$} algorithms at $95$\% confidence interval. As observed from Fig.~\ref{fig:nn}, the former two algorithms achieve significantly close solutions to {\sc OPT$_{\textsc{OPF}}$} for all case studies. The results illustrated in Fig.~\ref{fig:nn} show that the objective value of {\sc GRA} is nearly the same and sometimes slightly lower than that of {\sc GRA$_{\textsc{PF}}$} since {\sc GRA} employs a worst-case approach for quantifying the system losses. For optimizing the MG operation with $1400$ customers, {\sc GRA} required nearly $10$ milliseconds, whereas {\sc GRA$_{\textsc{PF}}$} terminated nearly in $20$ seconds. It is worthy to mention that the solutions obtained by Gurobi ({\sc OPT$_{\textsc{OPF}}$}) require approximately $60$ seconds for an MG with $1400$ customers. The running time of the candidate algorithms plotted in Fig.~\ref{fig:nnt} reflects the trade-off between optimality and scalability in the context of computational complexity.

\section{Conclusion}\label{sec:concl}

This paper presents a computationally efficient centralized algorithm for event-based DR management in MGs considering a large number of customers. The proposed algorithm relies on a greedy approach with a provable approximation guarantee. Various greedy strategies are analyzed and compared to identify the most efficient one in terms of quality of the solution. Furthermore, the proposed greedy approach is compared, in terms of computational time, to Gurobi optimizer known for its efficient computational speed. The results show that the Greedy Ratio Algorithm (GRA) is superior over the other strategies in terms of quality of the solution. The proposed algorithm is capable of optimizing the demand management in MGs with a large number of customers within acceptable time. For the majority of case studies the proposed algorithm demonstrated near optimal performance in little computational time which is an essential feature for event-based DR management in MGs. Lastly, we investigate the feasibility of the produced solutions by {\sc GRA} algorithm on a feeder from Canadian benchmark distribution system. The simulation results confirm that the produced solutions do not violate network constrains and voltage levels. Recently, this work has been extended to consider scheduling, optimal power flow and online setting \cite{KKEC16CSP, KCE2016OPF, instance1290}.

\medskip

\appendix 
In the appendix, the full proof of Theorem \ref{thm:alg-greedy} is presented followed by the formulation of {\sc UMDR-OPF} problem.

\subsection{Proof of Theorem \ref{thm:alg-greedy}}

\begin{proof} 
First, note that the problem is invariant, when the arguments of all complex-valued demands are rotated by the same angle (see Fig.~\ref{fig:rotate} for an illustration).

\begin{figure}[!ht]
	\begin{center} 
		\includegraphics[scale=0.8]{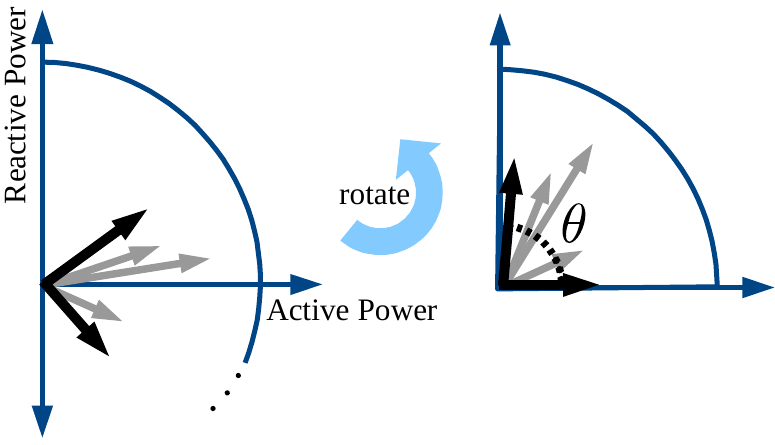}
	\end{center} 
	\caption{Each vector represents a power demand $S_k$. The figure shows that the demands are rotated by the same angle. $\theta$ is the maximum angle between any pair of demands. }
	\label{fig:rotate}
\end{figure}

\begin{figure*}[!t]
	\hrule
	\begin{align}
		&\frac{  \sum_{i=1}^n |d_i|^2 +  2\sum_{1\le i < j \le n} |d_i| \cdot |d_j|  } { \sum_{i =1}^n |d_i|^2 +  2\sum_{1\le i < j \le n} |d_i| \cdot |d_j| (\sin \theta_i \sin \theta_j + \cos \theta_i \cos \theta_j)} =\frac{ \sum_{i=1}^n |d_i|^2 +  2\sum_{1\le i < j \le n} |d_i| \cdot |d_j|  } { \sum_{i=1}^n |d_i|^2 +  2\sum_{1\le i < j \le n} |d_i| \cdot |d_j| \cos (\theta_i - \theta_j)} \label{eq:ind}\\
		&= \frac{ ( \sum_{i=1}^r |d_i|)^2 } { \displaystyle \sum_{i=1}^r |d_i|^2 + 2 \smashoperator{\sum_{1\le i<j<r}} |d_i|  |d_j| \cos (\theta_i - \theta_j) +  2 |d_r| \smashoperator{\sum_{1\le i<r}} |d_i|  \cos (\theta_i - \theta_r)}. \label{eq:den}
	\end{align}
	\hrule
\end{figure*}

Each customer's power demand has a phase angle between the reactive and apparent powers. Define $\theta$ to be the maximum difference between the phase angles of any pair of power demands of the customers (see Fig.~\ref{fig:rotate}).
Without loss of generality, this paper assumes that one of the demands, say $S_1$ is aligned along the positive real axis, and defines a class of sub-problems, by restricting the maximum phase angle $\theta$ (i.e., the argument) that any other demand makes with $S_1$ (see Fig.~\ref{fig:rotate} for an illustration). Note that in practice $\theta < \frac{\pi}{2}$, because there are regulations that require electric equipment to conform with a certain maximum power factor. More precisely, $\theta$ is usually restricted to be in the range of $[0, 36^{\circ}]$ \cite{1339347}. For the clarity of presentation, this paper assumes (via a rotation) that $P_k \ge 0$ and $Q_k \ge 0$.

\begin{figure}[!htb]
	\centering \vspace{-5pt}
	\includegraphics[scale=.8]{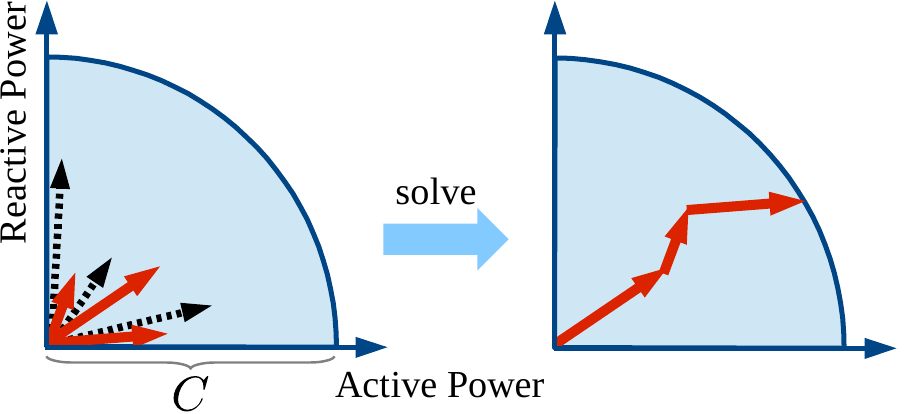} 
	\caption{The red vectors (thick arrows) represent a feasible solution to {\sc UMDR} such that the total magnitude of the red demands lies within the radius $C$.} 
	\label{fig:problem}
\end{figure}

Denote by $Z^{\rm GRA}$ the utility of the output solution of {\rm GRA}. Let $ S^\ast\subseteq {\cal N}$ be an optimal solution of {\sc UMDR}, and the $\Opt \triangleq \sum_{k \in  S^\ast} u_k$ be the corresponding total utility.
 
If  customers' complex-valued power demands are substituted in {\sc UMDR} problem by its real-valued magnitude and the binary decision variable $x_k$ is relaxed such that it takes non-negative real values instead of only integers values (i.e., $(x_k)_{k \in {\cal N}} \in [0, 1]^n$) the following linear programming (LP) problem is formulated as follows:

\begin{eqnarray}
\textsc{(LP1)} \qquad& \displaystyle \max_{x_k \in [0, 1]} \sum_{k\in{\cal N}}u_k x_k \label{LVmaxPA}\\
\text{subject to}\qquad & \displaystyle \sum_{k\in {\cal N}}|S_k| x_k \le C \label{C2}
\end{eqnarray}

Denote by $ S^\ast_{L} \subseteq {\cal N}$ an optimal solution of {\sc LP1} and by $\Opt_L \triangleq \sum_{k \in  S^\ast_{L}} u_k$ the corresponding total utility. In chapter 2 of \cite{KPP10book} it was shown that the optimal solution to {\sc LP1} problem can be determined easily, since the problem admits the greedy choice property. In other words, a global optimum of the problem can be achieved by choosing a series of locally optimal choices. The greedy choice for {\sc LP1} problem is to select customers in sorted order defined by Eqn.~\raf{eqn:gra}. Assume at some iteration $t$ adding the next customer to $X$ that causes capacity constraint violation, that is

\begin{equation}
\displaystyle \sum\limits_{k=1}^{t-1} |S_k| \le C \text{ and } \sum\limits_{k=1}^{t} |S_k| > C.
\end{equation}

The greedy execution is stopped at this point and the remaining capacity $C - \sum\limits_{k=1}^{t-1} |S_k|$ is occupied by the corresponding fractional part of the $t$-th customer's power demand. Observe that the preceding greedy strategy is the adapted version of {\sc GRA} algorithm for the relaxed {\sc UMDR} problem with continuous decision variable $(x_k)_{k \in {\cal N}} \in [0, 1]^n$ and complex-valued power demands. Let $ \hat{p} \triangleq \sum\limits_{k=1}^{t-1} u_k $ and $u_{\max} \triangleq \max_{k\in \cal N} u_k$. It was also shown in \cite{KPP10book} that

\begin{equation}
\label{eq:optl}
\displaystyle \Opt_L =  \hat{p} + (C - \sum\limits_{k=1}^{t-1} |S_k|)\frac{u_t}{|S_t|} \le \hat{p} + u_{\max}.
\end{equation}

Evidently, $Z^{\textsc {GRA}} \ge \hat{p}$. This gives

\begin{equation}
\label{eq:optl2}
\displaystyle \Opt_L \le Z^{\rm GRA} + u_{\max}.
\end{equation}

On the other hand, by Lemma \ref{lem:tb} below it follows that

\begin{equation}
\displaystyle \cos \tfrac{\theta}{2}\cdot \sum_{i \in S^\ast} |S_i| \le \Big| \sum_{i \in S^\ast} S_i \Big | \le C \, ,
\end{equation}

\noindent since  $\theta$ is restricted to be at most $\frac{\pi}{2}$. Note that the subset $S^{\ast}$ becomes a feasible solution to {\sc LP1} if  the relaxed decision variable is set $x_i=\cos \tfrac{\theta}{2}$ for $\forall i \in S^{\ast}$ and $x_i=0$ otherwise. This implies that

\begin{equation}
\label{eq:optl4}
\displaystyle \Opt_L \ge \cos \tfrac{\theta}{2} \cdot u(S^{\ast}) = \cos \tfrac{\theta}{2} \cdot \Opt.
\end{equation}

In {\rm GRA}, $Z^{\rm GRA} \ge u_{\max}$, and hence by Eqns~\raf{eq:optl2} and \raf{eq:optl4} it follows that

\begin{equation}
\label{eq:optl5}
\displaystyle Z^{\rm GRA}  \ge \frac{1}{2} \cos \tfrac{\theta}{2} \cdot \Opt.
\end{equation}
\end{proof}

\medskip

\begin{customlem}{1}
	\label{lem:tb}
	Given a set of 2D vectors $\{d_i \in \RR^2\}_{i=1}^n$
	$$ \frac{\sum_{i=1}^n |d_i| }{\big| \sum_{i =1}^n d_i \big|} \le \sec\tfrac{\theta}{2},$$
	where $\theta$ is the maximum angle between any pair of vectors  $\{d_i \in \RR^2\}_{i=1}^n$ and $0 \le \theta \le \frac{\pi}{2}$.
\end{customlem}
\begin{proof}
	The below statement will be proven by induction: $$\frac{(\sum_{i=1}^n |d_i| )^2}{|\sum_{i=1}^n d_i |^2} \le \frac{2}{\cos\theta + 1}$$  Note that $\sec \tfrac{\theta}{2} = \sqrt{\frac{2}{\cos \theta + 1}}$. 
	
	First,  expand the left hand side in Eqn.~\raf{eq:ind}, where $\theta_i$ is the angle that $d_i$ makes with the $x$ axis.

	Consider the base case: $n= 2$. Eqn.~\raf{eq:ind} is
	\begin{align}
	\frac{|d_1|^2 +|d_2|^2 + 2 |d_1|\cdot |d_2|}{|d_1|^2 +|d_2|^2 + 2 |d_1|\cdot |d_2| \cos(\theta)} = f\Big(\frac{|d_2|}{|d_1|}\Big),
	\end{align}
	where $$ f(x) \triangleq \frac{1+x^2+2x}{1+x^2+2x\cos \theta}.$$ The first derivative $f'(x)=$
	$$ \frac{(1+x^2+2x \cos \theta )(2x + 2) - (1+x^2 + 2x)(2x + 2 \cos \theta)}{(1 + x^2 + 2x \cos \theta )^2}$$
	 is zero only when $x=1$. Hence, $f(1)$ is an extreme point.  Compare $f(1)$ with $f(x)$ at the  boundaries $x\in \{0,\infty\}$: $f(1) = \frac{2}{\cos \theta + 1} \ge f(0) = \lim_{x \to \infty} f(x) = 1.$
	  Therefore, $f(x)$ has a global maximum of $\frac{2}{\cos \theta  + 1}$.
	
	Next proceed to the inductive step. Assume that $$\frac{\sum_{i=1}^{r-1} |d_i| }{\big| \sum_{i=1}^{r-1} d_i \big|} \le \sqrt{\frac{2}{\cos \theta + 1}},$$ where $r \in \{1,\ldots, n\}$. Without loss of generality, assume $\theta_1\ge\theta_3\ge \cdots \ge \theta_n \ge \theta_2.$
	
	Eqn.~\raf{eq:ind} can be rewritten to obtain Eqn.~\raf{eq:den}.


	Let $g(\theta_r)$ be the denominator of Eqn.~\raf{eq:den}. Take the  second derivative of $g(\theta_r)$:
	$$g''(\theta_r) = -2 |d_r| \sum_{1 \le i < r} |d_i| \cos(\theta_i - \theta_r).$$ Note that $\cos(\theta_i - \theta_r) \ge 0$, and hence, the second derivative is always negative. This indicates that all local extreme in $[0,\theta_{r-1}]$ of $g(\theta_n)$ are local maximum. Hence, the minimum occur at the boundaries:
	\begin{equation}
	\min_{\theta_r \in [0,\theta_{r-1}]}  g(\theta_r) \in \{g(0), g(\theta_{r-1})\}
	\end{equation}
	If $\theta_r \in \{0,\theta_r\}$ there must exist at least a pair of vectors in $\{d_i\}_{i=1}^r$ with the same angle. Combining these two vectors into one, one  can obtain an instance with $r-1$ vectors. Hence, by the inductive hypothesis, the same bound holds up to $r$ vectors. 
\end{proof}

\subsection{Formulation of {\sc UMDR-OPF} Problem}

The constrains present in OPF problem are broadly defined in the literature by two models, namely Bus Injection Model (BIM)  and Branch Flow Model (BFM) (a.k.a DistFlow model)~\cite{low2014convex1}. Variables in BIM (i.e., voltage and power) are assigned for every bus (or node). On the other hand, variables in BFM are assigned for every branch (or edge). BFM was first proposed by~\cite{baran1989placement}. It was shown by~\cite{equal} that both models are in fact equivalent. In this paper we utilize BFM model~\cite{low2014convex1} for introducing voltage, network and operational constrains into {\sc UMDR} problem.

In order to model the power flow constrains we represent an MG by a graph $\cG=(\cV,\cE)$. The set of nodes $\cV$ denote the buses, whereas the set of edges $\cE$ denote the distribution lines. We index the nodes in $\cV$ by $\{0, 1..., |\cV|\}$, where the node $0$ denotes the generation source or MG. For node $i \in \cV$, we denote its voltage by $V_i \in \CC$. For each edge $e=(i,j) \in \cE$, we denote its current from $i$ to $j$ by $I_{i,j}$, its transmitted power by $\hat{S}_{i,j}$, and its impedance by $z_{i,j} \in \CC$ (also denoted by $z_e$). Let $v_i \triangleq |V_i|^2$ and $\ell_{i,j} \triangleq |I_{i,j}|^2$ be the magnitude square of voltage and current respectively. For each node $i\in\cV\backslash\{0\}$, there is a set of customers attached to $i$, denoted by $\cN_i$. Given a complex number $\psi \in \CC$, denote its real and imaginary components by $\re(\psi)$ and $\im(\psi)$ respectively, and its complex conjugate by $\psi^\ast$. Let $v_{\min}, v_{\max} \in\RR^+$ be the minimum and maximum allowable voltage magnitude square at any node. We define the {\sc UMDR-OPF} problem by the following mixed integer non-linear programming problem.
	
\begin{align}
&\textsc{(UMDR-OPF)}\quad \max_{\substack{x_k, v_i,\ell_{i,j}, \hat{S}_{i,j} \;\;}} \sum_{k \in \cN} u_k x_k  \notag \\
\text{s.t.} \ & 	\ell_{i,j} \geq \frac{|\hat{S}_{i,j}|^2}{v_i}\, ,  \qquad \forall (i,j) \in \cE  \label{eq:c0}\\
&  \hat{S}_{i,j}=  \sum_{k \in \cN_j} S_k x_k  + \sum_{l: (j,l)\in \cE} \hat{S}_{j,l} + z_{i,j}\ell_{i,j}, \ \forall (i,j) \in \cE  \label{eq:S}\\
&	v_j = v_i + |z_{i,j}|^2 \ell_{i,j} - 2 \re(z_{i,j}^\ast  \hat{S}_{i,j}), \quad  \forall (i,j) \in \cE \label{eq:cv}\\
& |\hat{S}_{0,1}| \le C \label{eq:c1}\\
& v_{\min} \le v_j \le v_{\max},  \quad  \forall j \in\cV\backslash \{0\}\label{eq:c2}\\
& x_k \in \{0,1\},\qquad \forall k \in \cN\\
& v_i \in \RR^+, \ \ell_{i,j} \in \RR^+, \ \hat{S}_{i,j} \in \CC \label{clast}
\end{align}	

We remark that OPF problem is non-convex due to the quadratic equality constraints

\begin{eqnarray}
\ell_{i,j} = \frac{|\hat{S}_{i,j}|^2}{v_i}\, ,  \qquad \forall (i,j) \in \cE
\end{eqnarray}
and thus is computationally intractable in general. We therefore consider relaxing them to inequalities~\raf{eq:c0} in {\sc UMDR-OPF} problem to convexify the problem. Due to the same reason in~\cite{6980137,low2014convex1} the same relaxation is adopted. Obviously, when the equality in~\raf{eq:c0} is attained in the solution to {\sc UMDR-OPF} then the OPF relaxation is exact.

\bibliographystyle{IEEEtran}
\bibliography{reference}



\end{document}